\documentstyle{article}
\input amssym.def
\input amssym.tex

\author{Dmitri Akhiezer}
\title{\bf Real group orbits on flag manifolds }

\date{ }

\begin{document}
             
\maketitle

\abstract
\noindent
In this survey,
we gather together
various results on the action of a real form $G_0$ of a
complex semisimple group $G$ on its flag
manifolds. 
We start with the finiteness theorem of J.~Wolf
implying that at least one of the $G_0$-orbits is open.
We give a new proof of the converse statement for real forms of inner type,
essentially due to F.M.~Malyshev. 
Namely, if a real form of inner type $G_0 \subset G$ has an open orbit
on a complex
algebraic homogeneous space $G/H$ 
then $H$ is parabolic. In order to
prove this, we recall, partly with proofs,
some results of A.L.~Onishchik on the factorizations
of reductive groups. Finally, we discuss the cycle
spaces of open $G_0$-orbits and define
the crown of
a symmetric space of non-compact type. 
With some exceptions, the cycle space agrees
with the crown.
We sketch a complex analytic
proof of this result, 
due to G.~Fels, A.~Huckleberry and J.~Wolf.  

\section{Introduction}

The first systematic treatment
of the orbit structure of a complex flag manifold $X=G/P$
under the action of a real form $G_0 \subset G$
is due to J.~Wolf \cite{W}.
Forty years after his paper,
these real group orbits and their cycle spaces
are still an object of intensive research.
We present here some results in this area, together with
other related results
on transitive and locally transitive actions
of Lie groups on complex manifolds.

\medskip
\noindent
The paper is organized as follows.
In Section 2 we prove the celebrated finiteness theorem for $G_0$-orbits
on $X$ (Theorem 2.3). 
We also state a theorem
characterizing open $G_0$-orbits on $X$ (Theorem 2.4).
All results of Section 2 are taken from \cite{W}.
In Section 3 we recall for future use 
a theorem, due to B.~Weisfeiler \cite{We}
and A.~Borel and J.~Tits \cite{BT}. Namely, let
$H$ be an algebraic subgroup
of a connected reductive group $G$. Theorem 3.1 shows that
one can find
a parabolic subgroup 
$P \subset G$ containing $H$, such that the unipotent radical of  
$H$ is contained in the unipotent radical of $P$.
In Section 4 we consider the factorizations
of reductive groups. The results of this
section are due to A.L.~Onishchik \cite{On1, On2}. 
We take for granted his list
of factorizations $G = H_1\cdot H_2$,
where $G$ is a simple algebraic group over ${\Bbb C}$
and $H_1, H_2 \subset G$ are reductive complex subgroups
(Theorem 4.1), and deduce from it 
his theorem on real forms. Namely, a
real form $G_0$ acting locally transitively
on an affine homogeneous space $G/H$ 
is either $SO_{1,7}$ or $SO_{3,5}$.
Moreover, in that case $G/H =
SO_8/Spin_7$ and the action of $G_0$ is in fact transitive
(Corollary 4.7).
This very special homogeneous space 
of a complex group $G$ has on open orbit of 
a real form $G_0$,
the situation being typical for flag manifolds. 
One can ask what homogeneous spaces 
share this
property. It turns out that if a real form of inner type
$G_0 \subset G$ has an open orbit on a
homogeneous space $G/H$ with $H$ algebraic, then $H$
is in fact parabolic and so $G/H$ is a flag manifold.
We prove this in Section 5 (see Corollary 5.2)
and then retrieve the result of F.M.~Malyshev
of the same type
in which the isotropy subgroup is not necessarily algebraic
(Theorem 5.4). It should be noted that,
the other way around, the statement  
for algebraic homogeneous spaces can be deduced from
his theorem. Our proof of both results is new.

\medskip
\noindent
Let $K$ be the complexification
of a maximal compact subgroup $K_0 \subset G_0$.
In Section 6 we briefly recall the Matsuki correspondence 
between $G_0$- and $K$-orbits on a flag manifold.
In Section 7 we define,
following the paper of S.G.~Gindikin and the author \cite{AG}, the 
crown $\Xi $ of $G_0/K_0$ in $G/K$. We also introduce the cycle
space of an open $G_0$-orbit on $X=G/P$, first
considered by 
R.~Wells and J.~Wolf \cite{WeW},
and state a theorem describing the cycle spaces
in terms of the crown (Theorem 7.1). In fact, 
with some exceptions which are well-understood,
the cycle space of an open $G_0$-orbit on $X$
agrees with $\Xi $
and, therefore, is independent of the flag manifold. 
In Sections 8 and 9,
we give an outline of the original proof
due to G.~Fels, A.~Huckleberry and J.~Wolf \cite{FHW}, 
using the methods of complex analysis.
One ingredient of the proof is a theorem of
G.~Fels and A.~Huckleberry \cite{FH},
saying that $\Xi $ is a maximal $G_0$-invariant, Stein and
Kobayashi hyperbolic domain in $G/K$ (Theorem 8.4).
Another ingredient is the construction of the Schubert domain,
due to A.~Huckleberry and J.~Wolf \cite{HW} and
explained in Section 9.
Finally, in Section 10
we discuss complex geometric properties of flag domains.
Namely, let $q$ be the dimension
of the compact $K$-orbit in an open $G_0$-orbit.
We consider measurable open $G_0$-orbits and 
state the theorem of W.~Schmid and J.~Wolf \cite{SW}
on the $q$-completeness of such flag domains.

\medskip
\noindent
Given
a $K$-orbit $O$
and the corresponding $G_0$-orbit
$O^\prime $ on $X$,
S.G.~Gindikin and T.~Matsuki suggested
to consider the subset $C(O) \subset G$ of all $g\in G$, such that  
$gO\cap O^\prime \ne \emptyset$ and
$ gO \cap O^\prime $ is compact, see \cite{GM}.
If $O$ is compact then $O^\prime $ is open and $O \subset O^\prime $.
Furthermore,
in this case $C(O) = \{g\in G \ \vert \ gO \subset O^\prime\}$
is precisely the set whose connected component $C(O)^\circ$ at
$e \in G$ is the cycle space of $O^\prime $ lifted to $G$.
This gives a natural way of generalizing the notion
of a cycle space to lower-dimensional $G_0$-orbits.
Recently, using this generalization, T.~Matsuki
carried over Theorem 7.1 to
arbitrary $G_0$-orbits on flag manifolds, see \cite{M4} and Theorem 7.2.
His proof is beyond the scope of our survey. 
 
\medskip
\noindent

\section{Finiteness theorem}

Let $G$ be a connected 
complex semisimple Lie group, 
$\goth g$ the Lie algebra of $G$,
and $\goth g_0$ a real form of
$\goth g$. 
The complex conjugation of $\goth g$ over $\goth g_0$ is denoted
by $\tau $.
Let $G_0$ be the connected real Lie subgroup 
of $G$ with Lie algebra $\goth g_0$.
We are interested
in $G_0$-orbits on flag manifolds of $G$.
By definition, these manifolds are the quotients of the form $G/P$,
where $P \subset G$ is a parabolic subgroup.  
It is known that the intersection of two parabolic 
subgroups in $G$ contains a maximal torus of $G$. Equivalently,
the intersection of two parabolic subalgebras
in $\goth g$ contains a Cartan subalgebra of $\goth g$.
We want to prove a stronger statement in the case when
the parabolic subalgebras are $\tau $-conjugate.
We will use the notion of a Cartan subalgebra for an arbitrary 
(and not just semisimple) Lie algebra $\goth l$
over any field $k$.
Recall that a Lie subalgebra
$\goth j \subset \goth l$
is called a Cartan subalgebra if $\goth j$ is nilpotent
and equal to its own normalizer.
Given a field extension $k \subset k^\prime $,
it follows from that definition that
$\goth j$ is a Cartan subalgebra in $\goth l$
if and only if $\goth j\otimes _k k^\prime $ is a Cartan
subalgebra in $\goth l\otimes _k k^\prime $.
We start with a simple general observation. 

\medskip \noindent {\bf Lemma 2.1.}\ {\it
Let $\goth g$ be any complex Lie algebra,
${\goth g }_0 $ a real form of $\goth g$, and
$\tau : \goth g \to \goth g$ the complex conjugation of $\goth g$
over $\goth g_0$. 
Let $\goth h \subset \goth g$ be a complex Lie 
subalgebra. Then $\goth h \cap \goth g _0$ is a real form of 
$\goth h \cap \tau(\goth h)$.}

\medskip \noindent
{\it Proof.}\ For any $A \in \goth h \cap \tau (\goth h)$ one has
$2A = (A+\tau (A))+ (A - \tau (A))$, where the first summand is contained in
$\goth h \cap \goth g_0$ and the second one gets into that subspace after 
multiplication by $i$. 
\hfill{$\square $}

\medskip
\noindent
The
following corollary will be useful.

\medskip
\noindent
{\bf Corollary 2.2.}\ {\it If $\goth p$ is a parabolic subalgebra
of a semisimple algebra $\goth g$, then $\goth p \cap \tau (\goth p)$
contains a $\tau $-stable Cartan subalgebra $\goth t $ of $\goth g$.}

\medskip
\noindent
{\it Proof.}\ Choose a Cartan subalgebra   
$\goth j$ of
$\goth p \cap \goth g _0$. Its complexification $\goth t $
is a Cartan subalgebra
of $\goth p \cap \tau (\goth p)$, which is $\tau $-stable. 
Now, $\goth p \cap \tau (\goth p)$ contains a Cartan subalgebra 
$\goth t ^\prime $ of $\goth g$.
Since $\goth t$ and $\goth t^\prime $ are conjugate 
as Cartan subalgebras of $\goth p \cap \tau (\goth p)$, 
it follows that
$\goth t$ is itself 
a Cartan subalgebra of $\goth g$. 
\hfill {$\square $}

\medskip
\noindent
The number of conjugacy classes of Cartan subalgebras of a real semisimple 
Lie algebra is finite. This was proved independently by A.~Borel and B.~Kostant
in the fifties of the previous century, see \cite{Ko}. Somewhat later,
M.~Sugiura determined explicitly the number of conjugacy classes and
found their representatives for each simple Lie algebra, see \cite {Su}.
Let $\{\goth j_1, \ldots, \goth j_m\}$ be a complete system of representatives
of Cartan subalgebras of $\goth g_0$. For each $k, \ k=1,\ldots, m,$
the complexification $\goth t_ k $ of $\goth j_k$
is a Cartan subalgebra of $\goth g$. 

\medskip
\noindent
{\bf Theorem 2.3 ({\rm J.~Wolf \cite {W}, Thm.\,2.6}).}\ 
{\it For any parabolic subgroup $P \subset G$
the number of $G_0$-orbits on $X=G/P$ is finite.}

\medskip
\noindent{\it Proof.}\ Define 
a map $\iota : X \to \{1,\ldots, m\}$ as follows.
For any point $x \in X$ let $\goth p_x$ be the isotropy
subalgebra of $x$ in $\goth g$. By Corollary 2.2, we can choose a Cartan
subalgebra $\goth j _x$ of $\goth g _0$ in $\goth p _x$. 
Take $g\in G_0$ so that ${\rm Ad} g \cdot \goth j_x = \goth j_k$
for some $k, k=1, \ldots , m$. Since $\goth j_k $ and $\goth j_l$
are not conjugate for $k \ne l$, the number $k$ does not depend on
$g$. 
Let $k =\iota (x)$. Then $\iota (x)$ is constant along the orbit $G_0(x)$.
Now, for $\iota (x) $ fixed
there exists $g\in G_0$ such that $\goth p _{gx} $ 
contains $\goth t_k$ with fixed $k$. Recall that
a point of $X $ is uniquely determined by its isotropy subgroup.
Since there are only finitely many parabolic subgroups containing
a given maximal torus,
the fiber of $\iota $ has finitely many $G_0$-orbits.
\hfill {$\square $}

\medskip \noindent
As a consequence of Theorem 
2.3, we see that at least one $G_0$-orbit is open in 
$X$. We will need a description
of open orbits in terms of isotropy subalgebras of their points. 
Fix a Cartan
subalgebra $\goth t \subset \goth g$. Let $\Sigma = \Sigma (\goth g, \goth t)$
be the root system, 
$\goth g_\alpha \subset \goth g,\ \alpha \in \Sigma$, 
the root subspaces, $\Sigma ^+ = \Sigma^+(\goth g, \goth t)
\subset \Sigma $ a positive subsystem,
and $\Pi $ the set of simple roots
corresponding to $\Sigma ^+$. 
Every $\alpha \in \Sigma $ has a unique expression $\alpha = \sum 
_{\pi \in \Pi}\ n_{\pi}(\alpha )\cdot \pi $, where $n_{\pi}(\alpha )$
are integers, all non-negative for $\alpha \in \Sigma ^+  $
and all non-positive for $\alpha \in \Sigma ^- = -  \Sigma ^+$.
For an arbitrary subset $\Phi \subset \Pi $ 
we will use the notation
$$\Phi ^r =\{ \alpha \in \Sigma \ \vert \ n_{\pi }(\alpha ) = 0
\ {\rm whenever}\ \pi \not \in \Phi \},\ \ \Phi ^u
=\{ \alpha \in \Sigma ^+ \ \vert \ \alpha \not \in \Phi ^r\}.$$
Then the standard parabolic subalgebra $\goth p _\Phi\subset \goth g$
is defined by
$$\goth p _{\Phi } = \goth p _{\Phi }^r + \goth p_{\Phi }^u,$$
where
$$\goth p _{\Phi }^r = \goth t + \sum _{\alpha \in \Phi ^r}\ \goth g_{\alpha }
\ \
{\rm and}
\ \    
\goth p _{\Phi }^u = \sum _{\alpha \in \Phi ^u}\ \goth g _{\alpha }$$
are the standard reductive Levi subalgebra and, respectively, the 
unipotent radical of $\goth p_{\Phi }$. In the sequel,
we will also use the notation
$$\goth p_{\Phi}^{-u} = \sum _{- \alpha \in \Phi ^u} \goth g_{\alpha }.$$ 
Now, let 
$\goth k_0$ be a maximal compact subalgebra of $\goth g_0$.
Then we have the Cartan involution $\theta : \goth g_0 \to \goth g_0$
and the Cartan decomposition
$\goth g_0 = \goth k_0 + \goth m_0$,
where $\goth k_0$ and $\goth m_0$ are the eigenspaces of $\theta $ 
with eigenvalues 1 and, respectively, -1.
A $\theta $-stable Cartan subalgebra $\goth j \subset\goth g_0$
is called fundamental (or maximally compact) if $\goth j \cap \goth k_0$
is a Cartan subalgebra of $\goth k_0$. More generally,
a Cartan subalgebra $\goth j \subset \goth g_0$
is called fundamental if $\goth j$ is conjugate to a
$\theta $-stable fundamental Cartan subalgebra. 
It is known that any two fundamental Cartan subalgebras of $\goth g_0$
are conjugate under an inner automorphism of $\goth g_0$.
We will assume that a Cartan subalgebra $\goth t \subset \goth g$
is $\tau $-stable. In other words, $\goth t = \goth j ^{\Bbb C}$,
where $\goth j$ is a Cartan subalgebra in $\goth g_0$. Then $\tau $
acts on $\Sigma $ by
$\tau (\alpha ) (A) = \overline {\alpha (\tau\cdot A)}$,
where $\alpha \in \Sigma ,\, A \in \goth t$.

\medskip
\noindent
{\bf Theorem 2.4 ({\rm J.~Wolf \cite{W}, Thm.\,4.5}).}\
{\it Let $X = G/P$ be a flag manifold. Then the $G_0$-orbit
of $x_0 = e\cdot P$ is open in $X$ if and only if $\goth p = \goth p_{\Phi}$,
where

{\rm (i)} $\goth p \cap \goth g_0 $ contains a fundamental Cartan subalgebra 
$\goth j \subset \goth g_0$;

{\rm (ii)} $\Phi $ is a subset of simple roots for 
$\Sigma ^+ (\goth g, \goth t), 
\ \goth t = \goth j ^{\Bbb C}$,
such that $\tau \Sigma^+ = \Sigma ^-$.

}  
\medskip \noindent
The proof can be also found in \cite {FHW}, Sect. 4.2.

\section{Embedding a subgroup into a parabolic one }

Let $G$ be a group. The normalizer of a subgroup $H \subset G$
is denoted by $N_G(H)$. 
For an algebraic group $H$ the unipotent radical is denoted
by $R_u(H)$.

\medskip
\noindent
Let $U$ be an algebraic unipotent subgroup of a complex semisimple group $G$.
Set $N_1 = N_G(U)$, $U_1 = R_u(N_1) $, and continue inductively:

$$N_k = N_G(U_{k-1}),\ \ U_k = R_u(N_k),\ \ k \ge 2.$$
Then         
$U \subset U_1$ and
$ U_{k-1} \subset U_k ,\ \ N_{k-1}\subset N_k$
for all $k \ge 2$. Therefore there exists an integer $l$, such that
$U_l = U_{l+1}$. This means that $U_l$ coincides with
the unipotent radical of its normalizer. We now recall the
following general theorem of fundamental importance.

\medskip \noindent
{\bf Theorem 3.1 ({\rm B.~Weisfeiler \cite{We}, A.~Borel, J.~Tits \cite{BT},
Cor. 3.2}).}\ {\it Let $k$ be an arbitrary field, $G$ a connected
reductive
algebraic group defined over $k$, and $U$ a 
unipotent algebraic subgroup of $G$. 
If the unipotent radical of the normalizer
$N_G(U)$ coincides with $U$, then $N_G(U)$
is a parabolic subgroup of $G$.}  

\medskip \noindent 
For $k ={\Bbb C}$, which is the only case we need,
the result goes back to V.V.~Morozov, see \cite{BT}, Remarque 3.4.
In the above form, the theorem was conjectured by I.I.~Piatetski-Shapiro, see
\cite{We}. For future references, we state the following corollary
of Theorem 3.1.

\medskip \noindent
{\bf Corollary 3.2.}\ {\it Let $k = {\Bbb C}$ and
let $G$ be as above. The normalizer $N_G(U)$
of a unipotent algebraic subgroup $U \subset G$
embeds into a parabolic subgroup $P \subset G$ in such a way that
$U \subset R_u(P)$. For any algebraic subgroup $H \subset G$
there exists an embedding into a parabolic subgroup $P$,
such that $R_u(H) \subset R_u(P)$.}

\medskip \noindent
{\it Proof.} Put $P = N_G(U_l)$ in the above construction.
Then $U \subset U_l = R_u(P)$.
This proves the first assertion. To prove the second one,
it suffices to take $U = R_u(H)$.
\hfill$\square$

\section{Factorizations of reductive groups}

The results of this section are due to A.L.~Onishchik. 
Let $G$ be a group, $H_1, H_2 \subset G$ two subgroups.
A triple $(G;H_1,H_2)$ is called a factorization of $G$
if for any $g\in G$ there exist
$h_1 \in H_1$
and $h_2 \in H_2$, such that $g = h_1 \cdot h_2$. 
In the Lie group case a factorization
$(G;H_1,H_2)$ gives rise to the factorization $(\goth g; \goth h_1, \goth h_2)$
of 
the Lie algebra $\goth g$. By definition, this means that $\goth g =
\goth h_1 + \goth h_2$. Conversely, if $(\goth g; \goth h_1, \goth h_2)$
is a factorization of $\goth g$ then the product $H_1\cdot H_2$
is an open subset in $G$ containing the neutral element. In general,
this open set 
does not coincide with $G$, and so
a factorization $(\goth g; \goth h_1, \goth h_2)$ is sometimes called a local
factorization of $G$. But, if $G, H_1$ and $H_2$ are
connected reductive (complex or real) Lie groups then 
every local factorization is (induced by) a global one, 
see \cite {On2}. We will give a simple proof of this fact below,
see Prop. 4.3 and Prop. 4.4.

\medskip
\noindent
All factorizations of connected compact Lie groups are classified
in \cite {On1}, see also \cite {On3}, \S\,14. 
If $G, H_1$ and $H_2$ are connected reductive
(complex or real) Lie groups, then the same problem is solved in \cite{On2}.
The core of the classification is the complete list of factorizations
for simple compact Lie groups. We prefer
to state the result for simple algebraic groups over $\Bbb C$.
If both subgroups $H_1,H_2$  
are reductive algebraic, then the list
is the same as in the compact case.

\medskip \noindent
{\bf Theorem 4.1 ({\rm A.L.~Onishchik \cite {On1, On2}}).}
{\it If $G$ is a simple algebraic group over $k ={\Bbb C}$
and $H_1, H_2$ are proper reductive algebraic subgroups of $G$, then,
up to a local isomorphism and renumbering of factors, the factorization
$(G; H_1, H_2)$ is one of the following: 

{\rm (1)}\ $(SL_{2n}; Sp_{2n},
SL_{2n-1}), \ n\ge 2 $;

{\rm (2)}\ $(SL_{2n}; Sp_{2n},
S(GL_1\times GL_{2n-1}),\  n\ge 2$;

{\rm (3)}\ $(SO_7; G_2, SO_6)$;
 
{\rm (4)}\ $(SO_7; G_2, SO_5)$;

{\rm (5)}\ $SO_7; G_2, SO_3\times SO_2)$;

{\rm (6)}\ $(SO_{2n}; SO_{2n-1}, SL_n),\ n\ge 4 $

{\rm (7)}\ $(SO_{2n}; SO_{2n-1}, GL_n),\ n\ge 4$;

{\rm (8)}\ $(SO_{4n}; SO_{4n-1}, Sp_{2n}), \ n\ge 2$;

{\rm (9)}\ $(SO_{4n}; SO_{4n-1}, Sp_{2n}\cdot Sp_2),\ n\ge 2$;

{\rm (10)}\ $(SO_{4n}; SO_{4n-1}, Sp_{2n}\cdot k^*), \ n \ge 2$;

{\rm (11)}\ $(SO_{16}; SO_{15}, Spin _9)$;

{\rm (12)}\ $(SO_8; SO_7, Spin _7)$.}

\medskip
\noindent
Although this result is algebraic by its nature, the only known proof
uses topological methods. We want to show how Theorem 4.1 applies
to factorizations of complex Lie algebras involving their real forms.

\medskip
\noindent
{\bf Lemma 4.2.}\ {\it 
Let $\tau : \goth g \to \goth g$
be the complex
conjugation of a complex Lie algebra over its real form $\goth g_0$.
Let $\goth h \subset \goth g$ be a complex Lie subalgebra.
Then $\goth g = \goth g_0 + \goth h$ if and only if $\goth g = \goth h + 
\tau (\goth h)$.}

\medskip
\noindent
{\it Proof.}\ Let $\goth g = \goth g_0 + \goth h$.
For any $X \in {\goth g}_0$ one has $iX = Y + Z$,
where $Y \in \goth g_0$ and $Z \in \goth h$. This implies
$$2X=-iZ - \tau (iZ) \in \goth h +  \tau (\goth h).$$
Conversely, if $\goth g = \goth h + \tau (\goth h)$,      
then for any $X \in \goth g$ there exist $Z_1, Z_2 \in \goth h$, such that
$$X  = Z_1 + \tau (Z_2)
=(Z_1 - Z_2) + (Z_1 + \tau (Z_2)),$$
hence $X \in \goth h + \goth g_0$.
\hfill $\square $

\medskip
\noindent
{\bf Proposition 4.3.}\ {\it Let $G$ be a connected
reductive algebraic group over 
${\Bbb C}$ and let $H_1, H_2 \subset G$ be two reductive algebraic subgroups.
Then $\goth g = \goth h_1 + \goth h_2$ if and only if $G= H_1\cdot H_2$.}

\medskip
\noindent
{\it Proof.}\ It suffices to prove that the local factorization 
implies the global one. Let $X = G/H_2$ and let $n ={\rm dim}(X)$.
If $L $
is a maximal compact subgroup of $H_2$
and $K$ is a maximal compact subgroup of $G$, such that $L \subset K$,
then $X$ is diffeomorphic to a real vector bundle over $K/L$. 
Therefore $X$ is homotopically equivalent to a compact manifold
of (real) dimension $n$.                                              
On the other hand,
$H_1$ has an open orbit on $X$. Since $X$ is an affine variety,
closed $H_1$-orbits are separated by $H_1$-invariant regular functions.
But such functions are constant, so there is only one closed orbit.
Assume now that
$H_1$ is not transitive on $X$, so that the closed $H_1$-orbit
has dimension $m < n $. A well-known corollary of Luna's Slice Theorem 
displays $X$ as a vector bundle over the closed orbit, see \cite{Lu}. 
Thus $X$ is homotopically equivalent to that orbit and,
by the same argument as above, to a compact manifold
of (real) dimension $m$. 
Now, for a compact connected manifold $M$ of dimension $n$ one has
$H_i(M, {\Bbb Z}_2) = 0$ if $i > n$ and $H_n(M, {\Bbb Z}_2) 
\cong {\Bbb Z}_2 $,
see e.g. \cite{D}, Prop. 3.3 and Cor. 3.4.
Therefore two compact manifolds of dimensions $m$
and $n, m \ne n,$
are not homotopically equivalent,
and we get a contradiction.
\hfill $\square $

\medskip \noindent
As a corollary, we have a similar proposition for real groups.

\medskip \noindent
{\bf Proposition 4.4.}\ {\it Let $G, H_1$ and $H_2$ be real forms
of complex reductive algebraic 
groups $G^{\Bbb C}, H_1^{\Bbb C}$ and $H_2^{\Bbb C}$. 
For $G$ connected one has
$\goth g = \goth h_1 + \goth h_2$ if and only if $G=H_1\cdot H_2$.}

\medskip \noindent
{\it Proof.} 
If $\goth g = \goth h_1 + \goth h_2$ then 
$\goth g^{\Bbb C} = {\goth h_1}^{\Bbb C} + {\goth h_2}^{\Bbb C}$. Thus
$G^{\Bbb C} = H_1^{\Bbb C}\cdot H_2^{\Bbb C}$
by Proposition 4.3. The action of $H_1^{\Bbb C}
\times H_2^{\Bbb C}$ on $G^{\Bbb C}$, defined by
$$g \mapsto h_1gh_2^{-1},\ \ g \in G^{\Bbb C}, h_i \in H_i^{\Bbb C},$$
is transitive. For $g\in G \subset G^{\Bbb C}$ we have
the following estimate of the dimension of $(H_1 \times H_2)$-orbit
through $g$:
$${\rm dim}\, H_1gH_2 = {\rm dim}\, H_1 + {\rm dim}\, H_2 - {\rm dim} 
\, (H_1 \cap gH_2g^{-1}) \ge $$
$$ {\rm dim}_{\Bbb C}\, H_1^{\Bbb C}
+ {\rm dim}_{\Bbb C}\, H_2 ^{\Bbb C}
-{\rm dim} _{\Bbb C}\, H_1 ^{\Bbb C} \cap gH_2^{\Bbb C} g^{-1}
= {\rm dim}_{\Bbb C}\, G^{\Bbb C} = {\rm dim}\, G.$$
But $G$ is connected and each coset $H_1gH_2$ is open,
hence $G = H_1\cdot H_2$.
\hfill $\square $

\medskip
\noindent
We will use the notion of an algebraic subalgebra
of a complex Lie algebra $\goth g$, which corresponds
to an algebraic group $G$. A subalgebra $\goth h \subset \goth g$
is said to be algebraic, if the associated connected subgroup $H \subset G$
is algebraic. In general, this notion depends
on the choice of $G$. However, if $\goth g$ is semisimple,
which will be our case, then $\goth h$ is algebraic for some 
$G$ if and only if $\goth h$ is algebraic for any other $G$.
An algebraic subalgebra of $\goth g$ is said to be reductive,
if $H$ is a reductive algebraic subgroup of $G$. Again,
for $\goth g$ semisimple the choice of $G$ does not matter.

\medskip
\noindent 
{\bf Theorem 4.5 ({\rm cf. \cite{On2}, Thm. 4.2}).}
{\it Let $\goth g$ be a simple complex Lie algebra, $\goth h \subset \goth g$,
$\goth h \ne \goth g$,
a reductive algebraic subalgebra, and $\goth g _0$ a real form of $\goth g$.
If $\goth g = \goth g_0 + \goth h$ then $\goth g$
is of type ${\rm D}_4$, $\goth h$ is of type ${\rm B}_3$,
embedded as the spinor subalgebra, 
and $\goth g_0$ is either $\goth {so}_{1,7}$
or $\goth {so}_{3,5}$.} 

\medskip
\noindent
{\it Proof.}\ In the notation of Lemma 4.2 we have $\goth g = \goth h +
\goth \tau (\goth h)$. 
Note that $\tau (\goth h)$ is a reductive
algebraic subalgebra
of $\goth g$. 
Choose $G$ simply connected. Then $\tau $ lifts
to an antiholomorphic involution of $G$, which we again
denote by $\tau $. Let $H_1$ and $H_2$ be the connected
reductive algebraic subgroups of $G$ with Lie algebras $\goth h$ and,
respectively, $\tau (\goth h)$. By Proposition 4.3 we have
the global decomposition $G = H_1 \cdot H_2$. Since $H_1$ and $H_2$
are isomorphic, it follows from Theorem 4.1 that 
the factorization $(G; H_1, H_2)$ is obtained from factorization (12).
More precisely,
$G$ is isomorphic to
$Spin _8$, the universal covering group of $SO _8$, and $H_1, H_2$
are two copies of $Spin _7$ in $Spin _8$. We assume that $H_1$
is the image of the spinor representation $Spin _7 \to SO_8$
and $H_2$ comes from the embedding
$SO_7 \to SO_8 $.
The conjugation $\tau $ interchanges $H_1$ and $H_2$.
We want to replace $\tau $ by a holomorphic involutive automorphism
of $G$ with the same behaviour with respect to $H_1$ and $H_2$.
For this we need the following lemma.
   
\medskip
\noindent
{\bf Lemma 4.6.}\ {\it Let $G$ be a connected reductive algebraic group
over ${\Bbb C}$. Take a maximal
compact subgroup in $G$ which is invariant under $\tau $. 
Let 
$\theta : G \to G$ be the corresponding Cartan involution 
and let $\sigma = \theta \tau (= \tau \theta)$.
For a reductive algebraic subgroup $H \subset G$, the factorization
$(G; H, \tau (H))$ implies the factorization $(G; H, \sigma (H))$,
and vice versa.}

\medskip
\noindent  
{\it Proof.}\ 
First of all, if $(G; H_1, H_2)$ is a factorization
of a group 
then one also has the factorization
$(G; \tilde H_1, \tilde H_2)$, where 
$\tilde H_1 = g_1 H_1 g_1 ^{-1}, \tilde H_2 
= g_2 H_2 g_2 ^{-1}$
are conjugate subgroups.  
In the setting of the lemma, 
choose $\tilde H = gHg^{-1}$
so that a maximal compact subgroup of $\tilde H$ is contained in the
chosen maximal compact subgroup of $G$.
Then $\theta (\tilde H) = \tilde H$
and, consequently, 
$$\tau (H) \simeq \tau (\tilde H) = \sigma (\tilde H) \simeq \sigma (H),$$
where $\simeq $ denotes conjugation by an inner automorphism. 
By the above remark, one of the 
two factorizations
$(G;H , \tau (H)),\,(G; H, \sigma (H))$ implies the other.
\hfill $\square $

\medskip
\noindent
{\it End of proof of Theorem 4.5.}\ We can replace $H_2$
by a conjugate subgroup so that $H_1 $ and $H_2  $ are
interchanged by $\sigma $. The factorization is in fact defined for $SO_8$,
in which case the subgroups are only locally isomorphic. For
this reason, $\sigma $ is an outer automorphism. It follows that
the restriction of $\sigma $ to the real form,
i.e., the Cartan involution of the latter, is also an outer automorphism.
There are precisely two real forms of ${\rm D}_4$ with this property,
namely, ${\goth s \goth o}_{1,7}$ and ${\goth s \goth o}_{3,5}$.
The remaining non-compact real forms ${\goth s \goth o}_{2,6},
\ {\goth s \goth o}_{4,4}$, and ${\goth s \goth o}^* _8$
are of inner type, see Sect 5.
We still have to show that ${\goth s \goth o}_{1,7}$, as well as 
${\goth s \goth o}_{3,5}$, 
together with the complex spinor subalgebra,
gives a factorization of 
$\goth g = {\goth s \goth o}_8$. So let $\tau $ 
be the complex conjugation of $\goth g$
over ${\goth s \goth o}_{1,7}$ or ${\goth s \goth o}_{3,5}$. 
Define $\sigma $ as in the lemma 
and denote again by $\sigma $ the corresponding automorphism of $\goth g$.
The fixed point subalgebra of $\sigma $ has rank 3,
whereas $\goth g$ has rank 4. Thus $\sigma $
is an outer automorphism of $\goth g$.
There are three
conjugacy classes 
of subalgebras of type ${\rm B}_3$
in $\goth g$. Let $\Upsilon $ be the set of these conjugacy classes.
The group of outer isomorphisms
of $\goth g$ acts on $\Upsilon $ as the group
of all permutations of $\Upsilon $, isomorphic
to the symmetric group $S_3$. Choose ${\cal C }\in \Upsilon $ so that
$\sigma (\cal C) \ne \cal C$ and let $\goth h \in \cal C$.
Applying an outer automorphism of
$\goth g$, we can arrange that
$\goth h$ corresponds to $Spin _7$ and $\sigma (\goth h) = {\goth s \goth o}_7$.
Therefore $\goth g = \goth h + \sigma (\goth h)$ by Theorem 4.1
and $\goth g = \goth g_0 + \goth h$ by Lemma 4.6 and Lemma 4.2.
\hfill $\square $   

\medskip
\noindent
{\bf Corollary 4.7.}\ {\it Let $G$ be a simple algebraic group over ${\Bbb C}$,
$G_0$ a real form of $G$, and $H \subset G$ a proper 
reductive algebraic subgroup.
Then the following three conditions are equivalent:

{\rm (i)} $G_0$ is locally transitive on $G/H$;

{\rm (ii)} $G_0$ is transitive on $G/H$;

{\rm (iii)} up to a local isomorphism, $G = SO_8$, $H= Spin _7$,
${G_0} = SO_{1,7}$ or $SO_{3,5}$.}

\medskip
\noindent
{\it Proof.} Theorem 4.5 says that (i) and (iii) are equivalent.
Proposition 4.4 shows that (i) implies (ii).
\hfill $\square $

\section{Real forms of inner type}

Let $\goth g_0$ be a real semisimple Lie algebra
of non-compact type. Let
$\goth g_0 = \goth k_0 + \goth m_0$ be a Cartan
decomposition with the
corresponding Cartan involution $\theta $. It is known
that $\theta $ is an inner automorphism of $\goth g_0$ if and
only if $\goth k_0$ contains a Cartan subalgebra of $\goth g_0$.
If this is the case, we will say that the Lie algebra $\goth g_0$
and the corresponding Lie group $G_0$ is of inner type. 
Clearly, $\goth g_0  $
is of inner type if and only if all simple                   
ideals of $\goth g_0$ are of inner type. The Cartan classification
yields the following list of simple
Lie algebras
of inner type:
$${\goth s \goth l}_2(\Bbb R),
{\goth s \goth u}_{p,q},
{\goth s \goth o}_{p,q}\ (p\ {\rm or}\ q \ {\rm even}),
{\goth s \goth o}^*_{\, 2n},
{\goth s \goth p}_{\, 2n}({\Bbb R}),
{\goth s \goth p}_{p,q},$$
$${\rm EII, EIII, EV, EVI, EVII, EVIII, EIX, FI, FII, G}.$$
As we have seen in Sect.1, 
a conjugacy class of parabolic subalgebras
has a representative $\goth p$, such that $\goth g = \goth g_0 + \goth p$.
In other words, for any parabolic subgroup
$P \subset G$ the real form $G_0$ has an open orbit on $G/P$.
For real forms of inner type the converse is also true.

\medskip \noindent
{\bf Theorem 5.1.}\
{\it Let $\goth g$ be a complex semisimple Lie algebra, $\goth g_0$ 
a real form of $\goth g$ of inner type, and 
$\goth j$ a compact Cartan subalgebra of $\goth g_0$.
If $\goth h$ is an
algebraic Lie subalgebra of $\goth g$ satisfying 
$\goth g = \goth g_0 + \goth h$
then $\goth h$ is parabolic. Moreover, 
there exists an inner automorphism ${\rm Ad}(g),\ g\in G_0$,
such that $\goth h = {\rm Ad}(g)\cdot \goth p_{\Phi} $, where
$\Phi $ is a subset of simple roots for 
some ordering of
$\Sigma (\goth g, \goth j\, ^{\Bbb C})$. Conversely,
any such $\goth h$ satisfies $\goth g = \goth g_0 + \goth h$.}

\medskip \noindent
{\bf Corollary 5.2.}\ {\it Let $G$ be a complex semisimple
Lie group, $G_0\subset G$ a real form of inner type, 
and $H\subset G$ a complex algebraic subgroup.
If $G_0$ has an open orbit on $G/H$ then $H$ is parabolic.}

\medskip
\noindent For an algebraic Lie algebra $\goth h$
we denote by $R_u(\goth h)$ the unipotent radical
and by $L(\goth h)$ 
a reductive Levi subalgebra.
For the proof of Theorem 5.1 we will need a lemma which rules
out certain factorizations with semisimple factors.

\medskip
\noindent
{\bf Lemma 5.3.}\ {\it Let $\goth g$ be a simple complex Lie algebra
and let
$\goth h_1, \goth h_2 \subset \goth g$ be two semisimple real Lie subalgebras,
such that ${\goth h}_1 \cap {\goth h}_2 = \{0\}$.
Then $\goth g \ne {\goth h}_1 + {\goth h}_2$.}

\medskip
\noindent 
{\it Proof.}\ Assume ${\goth h}_1 +{\goth h}_2 = \goth g$.
Let $G$ be a simply connected Lie group with Lie algebra $\goth g$
and let $H_1, H_2$ be connected subgroups of $G$ with Lie algebras $\goth h_1$,
$\goth h_2$.
Then $G = H_1 \cdot H_2$ by Prop. 4.4. Therefore one can write 
$G$ as a homogeneous space $G = (H_1\times H_2)/(H_1\cap H_2)$,
where $H_1 \cap H_2$ embeds diagonally into the product.
Because $G$ is simply connected, we see that the
intersection $H_1 \cap H_2$ is in fact trivial.
But $H^3(G, {\Bbb R}) \cong {\Bbb R}$,
whereas  ${\rm dim}\, H^3(H_i, {\Bbb R}) \ge 1$,
see e.g. \cite{On3}, Ch.\,3, \S\,9, and
so the decomposition $G = H_1 \times H_2$
yields a contradiction.
\hfill $\square $ 

\medskip
\noindent
{\it Proof of Theorem 5.1.}\ 
Write $\goth g_0$ as the sum
of simple ideals  
$\goth g_{k,0},\, k=1,\ldots , m$. 
Each of them is stable under the Cartan involution $\theta $,
because $\theta $ is an inner automorphism. Furthermore,
each $\goth g_{k,0}$ is again of inner type. Thus the
complexification $\goth g _k= (\goth g_{k,0})^{\Bbb C}$ is a simple
ideal of $\goth g$,
and $\goth g = \goth g _1 \oplus \ldots 
\oplus \goth g _m$. Let $\pi _k  : \goth g \to \goth g _k$
be the projection map. 

Assume that $\goth h $ is reductive. We want to show that 
then $\goth h = \goth g$.
For each $k$ we have $\goth g _k = \goth g _{k,0} + \pi _k (\goth h)$.
If $\pi _k (\goth h) \ne \goth g _k$, then
$\goth g_{k,0}$ is isomorphic to ${\goth s \goth 0}_{1,7}$ or  
${\goth s \goth o}_{3,5}$ by Corollary 4.7. Since $\goth g_{k,0}$ is
of inner type, this can not happen. Hence $\pi _k (\goth h)
= \goth g _k$ for all $k$. In particular, $\goth h$ is semisimple, and
so we write
$\goth h$ as the sum of simple ideals $\goth h = \goth h_1
\oplus \ldots \oplus \goth h_n$. Since $\pi _k({\goth h}_l)$ is
an ideal in $\goth g_k$, there
are only two possibilities: $\pi _k ({\goth h}_l) = \goth g_k$
or $\pi _k ({\goth h}_l) =\{0\}$. If, for $k$ fixed, we have
$\pi _k (\goth h_l) \ne \{0\}$ and $\pi _k(\goth h _s) \ne \{0\}$
then in fact $ s= l$, because otherwise
$$ \goth g_k =[\goth g_k, \goth g_k] = [\pi _k(\goth h_l), \pi _k (\goth h_s)]
= \pi _k ([\goth h_l, \goth h_s]) =\{0\}.$$
We want to make sure that $m = n$. In that case, renumbering the 
simple ideals of $\goth g$, we get
$\goth h_l \subset \goth g_l$ for all $l$. This implies $\goth h_l =
\goth g_l$ for all $l$ and $\goth h = \goth g$.
Now, if $n < m $ then one and only one $\goth h_l$ 
projects isomorphically onto $\goth g_k$ and $\goth g _p$ for
$p \ne k$. Let $\omega _k =(\pi _k \vert \goth h_l)^{-1}$
and $\omega _p =(\pi _p \vert \goth h_l)^{-1}$. Then
$$\goth g_k \oplus \goth g _p = (\pi _k \oplus \pi _p)(\goth h_l) + 
(\goth g_{k,0} \oplus \goth g_{p,0}), $$
hence 
$$\goth h_l = \omega _k(\goth g_{k,0}) + \omega _p(\goth g_{p,0}),$$
and so a simple complex Lie algebra $\goth h_l$
is written as the sum of two real forms.
This contradicts Lemma 5.3.   

Assume from now on that $R_u(\goth h) \ne \{0\}$
and take an embedding of $\goth h$ into a parabolic subalgebra $\goth p$,
such that $R_u(\goth h) \subset R_u(\goth p)$, see Cor. 3.2.
Then $\goth g = {\goth g}_0 + \goth p$, i.e., the $G_0$-orbit
of the base point is open in $G/P$. By Theorem 2.4 $\goth p$
is a standard parabolic subalgebra, $\goth p = \goth p_{\Phi}$, where:

(i) $\goth p \cap {\goth g}_0 $ contains a fundamental Cartan subalgebra
$\goth j \subset {\goth g}_0$, which is now compact (recall that 
${\goth g}_0 $ is of inner type);

(ii) $\Phi $ is a subset of simple roots for some ordering of 
$\Sigma (\goth g, \goth t),\ \goth t = \goth j ^{\Bbb C}$ (since
$\goth j$ is compact, $\tau (\alpha) =
- \alpha $ for all $\alpha \in \Sigma (\goth g, \goth t)$
and $\tau \Sigma ^+ = \Sigma ^-$ for any choice of $\Sigma ^+$).

\smallskip \noindent By our construction,
$R_u(\goth h) \subset {\goth p}_{\Phi }^u$.
Applying an inner automorphism of $\goth p _{\Phi}$, 
assume that
$L(\goth h) \subset \goth p_{\Phi}^r$.
Since $\tau (\goth g_\alpha ) = \goth g _{-\alpha }$
for all root spaces, we have
$$\tau (\goth p_{\Phi}^r) = \goth p_{\Phi} ^r \ \ {\rm and}\ \  
\tau (\goth p_{\Phi}^u) = \goth p_{\Phi}^{-u}.$$
Observe that $\tau (\goth h) + \goth h = \goth g$ by Lemma 4.2. 
Therefore
$$\tau (R_u(\goth h)) + R_u(\goth h) = \tau (\goth p_{\Phi }^u) + 
\goth p_{\Phi }^u = \goth p _{\Phi }^{-u} + \goth p _{\Phi}^u,$$
and so we obtain
$$R_u(\goth h) = \goth p_{\Phi }^u.$$ 
We also have
$$\tau (L(\goth h)) + L(\goth h) = \goth p_{\Phi}^r,$$
hence, again by Lemma 4.2,
$$\goth p_{\Phi }^r = (\goth p _{\Phi}^r)_0 + L(\goth h), 
\ \ {\rm where} 
\ \ (\goth p _{\Phi}^r)_0 = \goth p_{\Phi }^r \cap \goth g_0.$$
Write $\goth p_{\Phi }^r = \goth s + \goth z$, where $\goth s$
is the semisimple part and $\goth z$ the center of $\goth p_{\Phi }^r$,
denote by $\pi _s, \pi _z$ the corresponding projections, and put
$\goth s _0 = \goth s \cap \goth g _0,\ \ \goth z_0 = 
\goth z \cap \goth g _0$. Then             
$\goth s = \goth s _0 + \pi _s(L(\goth h))$.
Since $\goth s _0$ is a semisimple algebra of inner type and
$\pi _s(L(\goth h))$ is reductive, we get $\pi _s(L(\goth h)) = \goth s$
by what we have already proved. Therefore $L(\goth h) = \goth s +\goth z_*$,
where $\goth z_* $ is an algebraic subalgebra in $\goth z$.
On the other hand, $\goth z = \goth z_0 + \pi _z(L(\goth h))
= \goth z_0 + \goth z_*$.  
But $\goth z_0$ is compact, so $\goth z_* = \goth z$
and $L(\goth h) = \goth p_{\Phi }^r$. Together with the
equality $R_u(\goth h )=
\goth p _{\Phi }^u$, this gives $\goth h = \goth p_{\Phi}$.
To finish the proof, recall that any two compact
Cartan subalgebras of $\goth g_0$ are conjugate by an inner automorphism.
For the converse statement of the theorem, 
note that, $\goth j$ being
compact, 
(ii) in Theorem 2.4 is fulfilled for any ordering of $\Sigma (\goth g,
\goth t)$.  
\hfill $\square $

\medskip
\noindent
We now recover a theorem of F.M.~Malyshev,
in which $\goth h$ is not necessarily algebraic.
Of course, our Theorem 5.1 is a special case of his result.
We want to show
that  
the general case can be obtained from that special one. 
We adopt the notation introduced in the above proof.
Namely,                                                      
$\goth s = \goth s_{\Phi }$ and $\goth z = \goth z_{\Phi}$ are the semisimple
part and, respectively,
the center of the reductive algebra $\goth p_{\Phi}^r$.

\medskip
\noindent 
{\bf Theorem 5.4\ ({\rm F.M.~Malyshev \cite {Ma1}}).}\  
{\it Let $\goth g$, $\goth g_0$ 
and 
$\goth j$ be as in Theorem 5.1.
If $\goth h$ is a complex
Lie subalgebra of $\goth g$ satisfying 
$\goth g = \goth g_0 + \goth h$
then  
there exists an inner automorphism ${\rm Ad}(g),\ g\in G_0$,
such that $\goth h = 
{\rm Ad}(g) (\goth a + \goth s _{\Phi}+ \goth p_{\Phi}^u)$, where
$\Phi $ is a subset of simple roots for 
some ordering of
$\Sigma (\goth g, \goth j\, ^{\Bbb C})$ and
$\goth a$ is a complex subspace of $\goth z _{\Phi}$ 
which projects onto the real form $(\goth z_{\Phi})_0$. 
Conversely,
any such $\goth h$ satisfies $\goth g = \goth g_0 + \goth h$.}

\medskip \noindent
{\it Proof.}\ Let $\goth h_{\rm alg}$ be the algebraic closure of $\goth h$,
i.e., the smallest algebraic subalgebra of $\goth g$ containing $\goth h$.
According to a theorem of C.~Chevalley, the commutator algebras
of $\goth h$ and $\goth h_{\rm alg}$ are the same, see \cite{C}. Chap.\,II,
Th\'eor\`eme 13. Applying
an inner automorphism ${\rm Ad}(g),\ g\in G_0$,
we get
$$\goth h_{\rm alg} = \goth p_{\Phi}=
\goth z_{\Phi} + \goth s_{\Phi} + \goth p_{\Phi}^u$$
by Theorem 5.1. Since $\goth h$
contains
$[\goth h_{\rm alg}, \goth h_{\rm alg}] = \goth s_{\Phi} + \goth p_{\Phi }^u$,
it follows that
$$\goth h = \goth a + \goth s_{\Phi} + \goth p_{\Phi}^u,$$
where $\goth a \subset \goth z_{\Phi}$ is a complex subspace. 
Observe that
$$\tau (\goth h) = \tau (\goth a) + \goth s_{\Phi} + \goth p_{\Phi}^{-u}.$$
By Lemma 4.2 we have
$\goth g = \goth h + \tau (\goth h)$.
Clearly, $\goth z_{\Phi}$ is $\tau $-stable. 
The above expression for $\tau (\goth h)$ shows that
$\goth z_{\Phi} = \goth a + \tau (\goth a)$.
Again by Lemma 4.2,
this implies $\goth z_{\Phi} = (\goth z_{\Phi})_0 + \goth a$ or,
equivalently,
$\goth z_{\Phi} = i\cdot (\goth z_{\Phi})_0 + \goth a$.
Thus $\goth a$ projects onto $(\goth z_{\Phi})_0$.
Since $\tau \goth p_{\Phi}^u = \goth p_{\Phi}^{-u}$,
the converse statement is obvious.
\hfill $\square $
 
\medskip \noindent
If $\goth g_0$ is a real form of outer type (= not 
of inner type), then a Lie subalgebra
$\goth h \subset \goth g$, satisfying $\goth g = \goth g_0 + \goth h$,
is in general very far from being parabolic. Some classification of such 
$\goth h$ is known for type ${\rm D}_n$, see \cite{Ma2}. Here is
a typical example of what can happen for other Lie algebras.

\medskip \noindent
{\bf Example 5.5.}\ Let $\goth g =\goth s \goth l _{2n}(\Bbb C),\ n> 1,$ and
let $\tau (A) = \bar A$ for $A \in \goth g$, so that $\goth g_0 = 
\goth s \goth l _{2n}(\Bbb R)$.
Then there is a fundamental Cartan subalgebra $\goth j \subset \goth g_0$
and an ordering of the root system $\Sigma (\goth g, \goth t),\ \goth t =
\goth j \, ^{\Bbb C},$ such that the set of simple roots 
$\Pi$ is of the form $\Pi = \Phi \sqcup \Psi \sqcup \{\gamma \}$,
where $\Phi $ and $\Psi $ are orthogonal,
$\tau(\Phi ) = -\Psi$ and $\tau (\gamma ) = -\gamma $.
The standard Levi subalgebra of $\goth p_{\Phi \sqcup \Psi}$
can be written as
$$\goth p_{\Phi \sqcup \Psi}^r = \goth s_1+ \goth s_2 + \goth z,$$ 
where $\goth s_1 $ and $\goth s_2$ are isomorphic simple algebras of type 
${\rm A}_{n-1}$
interchanged by $\tau $ and $\goth z$ is a $\tau$-stable one-dimensional torus.
Set
$$\goth h = \goth s_1 + \goth z + \goth p_{\Phi \sqcup \Psi}^u,$$
then
$$\tau (\goth h) = \goth s_2 + \goth z + \goth p_{\Phi \sqcup \Psi} ^{-u}.$$
Therefore $\goth h + \tau (\goth h) = \goth g$,
showing that $\goth g = \goth g_0 + \goth h$. Note that $\goth h$ is an
ideal in the parabolic subalgebra $\goth p = {\goth p}_{\Phi \sqcup \Psi}$,
such that $\goth p/\goth h$ is a simple algebra.

The construction of $\goth j$ and the ordering in $\Sigma (\goth g, \goth t)$
goes as follows. Take the Cartan decomposition 
$\goth g_0 = \goth k_0 + \goth m_0$, where $\goth k_0 = 
\goth s \goth 0_{2n}(\Bbb R)$. Define $\goth j$ as 
the space of block matrices                                                      

$$        \left(
           \matrix{ a_1  & b_1  & &  & & & \cr
                  -b_1  & a_1  & &  & & & \cr
              & & a_2 & b_2 & & & \cr
              & & -b_2 & a_2 & &  &\cr  
              &     & & &  \ddots & &\cr
              &     & & & & a_n & b_n \cr
              &     & & & & -b_n & a_n \cr
              } \right)$$

\noindent with real entries and
$\Sigma \,a_i = 0 $. Then $\goth j$ is a fundamental Cartan subalgebra and                                                       
$\goth j = \goth j \cap \goth k_0 + \goth j \cap \goth m_0$. Consider
$a_i$ and $b_i$ as linear functions on $\goth j$ and $\goth t =
\goth j\, ^{\Bbb C}$. Then it is easy
to determine the root system $\Sigma (\goth g, \goth t)$.
We list the roots that we declare positive:
$$i(b_p - b_q) \pm (a_p - a_q),\ i(b_p + b_q) \pm (a_p - a_q)\ \ (p< q),
\ {\rm and } \ 2ib_p\ \ (p,q=1,\ldots,n).$$
Let
$\Phi = \{\alpha_1, \ldots , \alpha_{n-1}\},
\Psi = \{\beta _1, \ldots , \beta _{n-1}\}$, where
$$\alpha _p = i(b_p - b_{p+1}) + a_p - a_{p+1},\
  \beta _p = i(b_p - b_{p+1}) - a_p + a_{p+1}\ (p=1,\ldots, n-1),$$
and let
   $\gamma = 2ib_n$.
Then the set of simple roots $\Pi$ is the union $\Pi = \Phi \sqcup \Psi
\sqcup \{\gamma \}$, $\Phi $ and $\Psi $ are orthogonal,
$\tau (\alpha _p) = - \beta _p$ for all $p$ and $\tau (\gamma) = -\gamma$.

\section{Matsuki correspondence}
Recall that $G_0$ is a real form of
a complex semisimple group $G$ and both $G_0$ and $G$ are connected.
Let $\goth g_0= \goth k_0 + \goth m_0$ be a Cartan decomposition,
$\goth k$ the complexification in $\goth g$ of $\goth k_0$, and $K$
the corresponding connected reductive subgroup of $G$. 

\bigskip
\noindent
{\bf Theorem 6.1 {(\rm T.~Matsuki \cite{M2})}.}
{\it
Let $O$ be a $K$-orbit and let $O^\prime $ be a $G_0$-orbit on $G/P$,
where $P \subset G$ is a parabolic subgroup.
The relation
$$O \leftrightarrow O^\prime \ \ \iff \ 
\ O\cap O^\prime \ne \emptyset {\rm \ \ and } \ O\cap O^\prime 
{\rm \ is\ compact} $$
defines a bijection between $K\setminus G\,/P$ and $G_0\setminus G\,/P$.}

\medskip
\noindent
A geometric proof of this result, using the moment map technique,
is found in \cite{MUV}, \cite{BL}. Note that $K$ is a spherical subgroup
of $G$, 
i.e., a Borel subgroup of $G$ has an open orbit on $G/K$.
It that case $B$ has finitely many orbits on $G/K$,
see \cite{Br}, \cite{Vin}.
Thus the set $K\setminus G\,/P$ is finite,
and so $G_0\setminus G\,/P$ is also finite
(another proof of Theorem 2.3). 

\medskip
\noindent
It can happen that both $K\setminus G\,/P$ and $G_0\setminus G\,/P$
are one-point sets. For $G$ simple, there are only two types
of such actions.

\medskip
\noindent
{\bf Theorem 6.2 {(\rm A.L.~Onishchik \cite{On2}, Thm. 6.1)}.} 
{\it If $G$ is simple and $G_0$ or, equivalently,
$K$ is transitive on $X=G/P$ then,
up to a local isomorphism,
                                                      
{\rm (1)}\ $G= SL_{2n}(\Bbb C),\ K= Sp_{2n}(\Bbb C),\ G_0 = SU^*_{2n},
\ X = {\Bbb P}^{2n-1}(\Bbb C)$, or 

{\rm (2)}\ $G = SO_{2n} (\Bbb C),\ K = SO_{2n-1}(\Bbb C),\ G_0 = SO^o_{2n-1,1},
\ X= SO_{2n}(\Bbb R)/U_n.$                            }

\medskip
\noindent
There are two important cases of the correspondence 
$O \leftrightarrow O^\prime $, namely, when one of the two orbits is open
or when it is compact. The first of the following two propositions is evident,
and the second one is due to T.~Matsuki \cite{M1}. 

\medskip
\noindent
{\bf Proposition 6.3.}\ {\it If $O$ is open then $O^\prime $ is compact
and $O^\prime \subset O$. If $O^\prime $ is open then $O$ is compact
and $O \subset O^\prime $.}
\hfill $\square $

\medskip \noindent
{\bf Proposition 6.4.}\ {\it If $O$ is compact then $O^\prime $ is open
and $ O \subset O^\prime $. If $O^\prime $ is compact then $O$ is open
and $O^\prime \subset O$.}

\medskip
\noindent
{\it Proof.}\ We prove the second statement. The proof of the 
first one is similar. Take a base point $x_0 \in O \cap O^\prime $
and let $P$ be the isotropy subgroup of $x_0$. 
Note that $G_0\cap P$ 
has only finitely many connected components,
since it is an open subgroup of a real algebraic group. By a
theorem of D.~Montgomery \cite{Mo}, $K_0$
is transitive on the compact homogeneous space $G_0/(G_0\cap P)$, hence
$$\goth g_0 = \goth k_0 + \goth g_0 \cap \goth p \subset \goth k + \goth p.$$ 
On the other hand, $\goth k_0 + i\goth m_0$
is the Lie algebra of a maximal compact subgroup of $G$,
which is transitive on $G/P$. Therefore
$$\goth g = \goth k_0+i\goth m_0 + \goth p$$
or, equivalently
$$\goth g = i\goth k_0 + \goth m_0 + \goth p,$$
and it follows that
$\goth g \subset \goth k + \goth g_0 + \goth p \subset \goth k + \goth p$,
i.e., $\goth g = \goth k+ \goth p$.
This means that $O = K(x_0)$ is open in $G/P$,
and the inclusion $O^\prime \subset O$ follows from Proposition 6.3.
\hfill $\square $

\section{Cycle spaces}

First, we recall the definition of the complex crown
of a real symmetric space $G_0/K_0$, see \cite {AG}.
Let $\goth a \subset \goth m_0$ be a maximal abelian subspace and let
$\goth a^+ \subset \goth a$ be the subset given by the inequalities $\vert 
\alpha (Y)
\vert < {\pi \over 2}$, where $Y \in \goth a$ and $\alpha $ runs
over all restricted roots,
i. e., the roots of $\goth g_0$ with respect to $\goth a$. Then the
crown is the set
$$\Xi = G_0 ({\rm exp} \, i\goth a^+) \, o \subset G/K,$$
where $o = e\cdot K \in G/K $ is the base point. 
The set $\Xi $ is open and the $G_0$-action on $\Xi $
is proper, see \cite {AG}. We discuss some 
properties of the complex manifold $\Xi $ in the next section.
Because
all maximal abelian subspaces in $\goth m_0$ are $K_0$-conjugate,
it follows that
$\Xi $ is independent of the choice of $\goth a$ and is therefore
determined by $G_0/K_0$ itself. Some authors call $\Xi$
the universal domain, see \cite{FHW}. We reserve this term
for the lift of $\Xi $ to $G$ and define the universal domain by
$$\Omega = G_0 ({\rm exp} \, i\goth a^+) K \subset G, $$
due to the properties of $\Omega $ which will soon
become clear. Of course, $\Omega $ is invariant under the right $K$-action
and $\Omega /K = \Xi $.

\medskip
\noindent
Next, we define the (linear) cycle space for an open $G_0$-orbit on $X=G/P$,
see \cite{WeW}.
Since full cycle spaces (in the sense of D.Barlet) are not discussed here,
we will omit the adjective ``linear".
Let $D$ be such 
an orbit and let $C_0$ be the corresponding
$K$-orbit, so that if $O\leftrightarrow O^\prime $
for $O = C_0$ and $O^\prime = D$. The orbit $C_0$ is a compact complex
manifold contained in $D$. 
Consider the open set $$G\{D\} = \{g\in G \ \vert \ gC_0 \subset D\}
\subset G$$              
and denote by $G\{D\}^\circ $ its connected component containing $e \in G$.
Observe that $G\{D\}$ is invariant under the right
multiplication by $L = \{ g\in G \ \vert \ gC_0 = C_0\}$
and left multiplication by $G_0$.
Since $L$ is 
a closed complex subgroup of $G$,
we have a natural complex structure on $G/L$. 
By definition, the cycle space 
${\cal M}_D$ of $D$ is the connected
component of $C_0(=e\cdot L)$ in
$G\{D\}/L$
with the inherited $G_0$-invariant complex structure.

\medskip
\noindent
In what follows, we assume $\goth g$ simple.
We will say that $G_0$ is
of Hermitian type if the symmetric space $G_0/K_0$
is Hermitian.
If this is the case then $\goth g$ 
has three irreducible components as $({\rm ad}\, {\goth k})$-module, 
namely,
$\goth g = \goth s^- + \goth k + \goth s^+$,
where $\goth s ^+, \goth s^-$
are abelian subalgebras. The subalgebras $\goth k + \goth s^+$
and $\goth k + \goth s^-$ are
parabolic. The corresponding parabolic subgroups are denoted $P^+$
and $P^-$. 
We have two flag manifolds
$X^+ = G/P^+ ,\ X^-=G/P^-$ with base points $x^+ = e\cdot P^+,\  x^-=
e\cdot P^-$
and two $G_0$-invariant complex structures on $G_0/K_0$
defined by the equivariant embeddings 
$g\cdot K_0 \mapsto gx^\pm \in X^\pm $. 
Each of the two orbits ${\cal B} = G_0(x^+)$ and $\bar {\cal B} = G_0(x^-)$  
is biholomorphically isomorphic
to the bounded symmetric domain associated to $G_0$.
The Lie algebra $\goth l$ of $L$ contains $\goth k$. 
If $G_0$ is of Hermitian type
and $\goth l$ coincides with $\goth p^+$
or $\goth p^-$ then we say that $D$ and, also,
the corresponding compact $K$-orbit $C_0$ is of (Hermitian)
holomorphic type.
If $G_0$ is of non-Hermitian type
then $\goth k$
is a maximal proper subalgebra of $\goth g$. Thus,
if $\goth l \ne \goth g$ then 
$\goth l = \goth k$. 
For $G_0$ of Hermitian type, each
flag manifold has exactly two $K$-orbits of holomorphic type.
All other $K$-orbits for $G_0$ of Hermitian type
and all $K$-orbits for $G_0$ of non-Hermitian type
are said to be of non-holomorphic type.
In the following theorem, we exclude the actions
listed in
Theorem 6.2. The symbol $\simeq $
means a $G_0$-equivariant biholomorphic isomorphism.

\medskip
\noindent
{\bf Theorem 7.1 {(\rm G.~Fels, A.~Huckleberry, J.A.~Wolf \cite{FHW}, 
Thm. 11.3.7)}.} 
{\it \ Assume $G$ simple and suppose $G_0$
is not transitive on $X =G/P$.
Let $D$ be an open $G_0$-orbit on $X$. 
If $D$ is of holomorphic type 
then ${\cal M}_D \simeq {\cal B}$ or ${\cal M}_D \simeq \bar
{\cal B}$.
In all other cases, 
$G\{D\}^\circ$ coincides with
the universal domain $\Omega \subset G$. Moreover, 
$\pi: G/K \to G/L$ is a finite covering map, which induces
a $G_0$-equivarint biholomorphic map $\pi\vert _{\Xi}: \Xi \to {\cal M}_D$.}

\medskip
\noindent
If $G_0$ is of Hermitian type then $\Xi \simeq {\cal B}\times \bar
{\cal B}$, see \cite{BHH}, Sect.3, \cite{GM},
Prop. 2.2, or \cite{FHW}, Prop. 6.1.9. The cycle space in that case
was first described by J.~Wolf and R.~Zierau \cite{WZ1,WZ2}. 
Namely,
in accordance with the above theorem, 
${\cal M}_D \simeq {\cal B}\times \bar {\cal B}$  
if $D$ is of non-holomorphic type and ${\cal M}_D \simeq {\cal B}$
or ${\cal M}_D \simeq \bar {\cal B}$ if $D$ is of holomorphic type.
For $G_0$ of non-Hermitian type,
the crucial equality $G\{D\}^\circ = \Omega $
is proved by G.~Fels and A.T.~Huckleberry in \cite{FH},
see Thm. 4.2.5,
using Kobayashi hyperbolicity of certain $G_0$-invariant
domains in $G/K$. 
In the next section, we consider
some properties of the crown $\Xi$,
which are important for that proof and are interesting in themselves.
After that,
we explain the strategy of their proof, without going into the details. 

\medskip
\noindent
Meanwhile,
the notion of the cycle space has been generalized to lower-dimensional
orbits and it turned out that its description in terms of
the universal domain holds in this greater generality. Namely,
given any $K$-orbit $O$ on $X=G/P$,
S.G.~Gindikin and T.~Matsuki \cite{GM}
defined a subset of $G$ by
$$C(O) = \{g \in G \ \vert \ gO \cap O^\prime \ne \emptyset {\rm \ and \ }
gO \cap O^\prime {\rm \ is \ compact}\},$$ 
where $O^\prime $ is the corresponding $G_0$-orbit, i.e.,
$O \leftrightarrow O^\prime $. 
Let $C(O)^\circ $ be the connected component of $C(O)$ containing $e \in G$.
Of course, if $D = O^\prime $
is open, then $C(O) = G\{D\}$ is
the open set considered above. 
The following theorem was stated as a conjecture
in \cite{GM}, see Conjecture 1.6. 

\medskip
\noindent
{\bf Theorem 7.2 {(\rm T.~Matsuki \cite{M4})}.}\ {\it Let $G$, $G_0$ and $X$ 
be as above. 
Then $C(O)^\circ = \Omega $ 
for all $K$-orbits on $X$ of non-holomorphic
type.}

\medskip
\noindent
{\it Remark.} The proof in \cite{M4}
uses combinatorial description of the
inclusion relations between the closures of $K$-orbits
on the flag manifolds of $G$. As a corollary, we get that $C(O)^\circ $
is an open set, which is not clear a priori.
If this is known,  
then Theorem 7.2 follows from \cite{HN} or
from Theorem 12.1.3 in \cite{FHW}. The latter
asserts that the connected component of 
the interior of $C(O)$, containing the neutral element $e\in G$,
coincides with $\Omega $.

\section{Complex geometric properties of the crown }

The following theorem proves the conjecture
stated in \cite{AG}.

\medskip
\noindent
{\bf Theorem 8.1 {(\rm D.~Burns, S.~Halverscheid, R.~Hind \cite{BHH})}.}
\ {\it The crown $\Xi $ is a Stein manifold.
}

\medskip 
\noindent
The crucial ingredient of the proof is the
construction of a smooth strictly plurisubharmonic function on $\Xi $          
that is $G_0$-invariant and gives an exhaustion of the orbit space 
$G_0\setminus \Xi$. We call such a function a BHH-function.
Let $\Gamma \subset
G_0$ be a discrete cocompact subgroup acting freely on $G_0/K_0$.
Then
$\Gamma $ acts properly and freely on $\Xi $ and
any BHH-function induces a plurisubharmonic exhaustion of $\Gamma 
\setminus \Xi$. Thus $\Gamma \setminus \Xi$ is a Stein manifold
and its covering $\Xi $
is also Stein.

\medskip
\noindent
We now want to give another application of BHH-functions.
Let $G_0 = K_0A_0N_0$ be an Iwasawa decomposition and
let $B$ be a Borel subgroup of $G$ containing 
the solvable subgroup $A_0N_0$. 
Then $B$ is called an Iwasawa-Borel
subgroup, the orbit $B(o) \subset G/K$
is Zariski
open and its complement, to be denoted
by ${\cal H}$, is a hypersurface. 
The set
$$\Psi = \bigcap _{g\in G_0} \, gB(o) = \bigcap _{k\in K_0}\, kB(o)$$
is open as the intersection of a compact family of open sets.
Let $\Xi_I $
be the connected component of $\Psi $ containing $o$.
L.~Barchini \cite{Ba} showed that 
$\Xi_I \subset \Xi$. The reverse inclusion was checked
in many special cases including all classical groups and all real forms
of Hermitian type, see
\cite {GM,KS}. The proof in the general case
is due to A.~Huckleberry, see \cite{Ha, FH} and \cite{FHW},
Remark 7.2.5. His argument is as follows.
It is enough to prove that ${\cal H}\cap \Xi = \emptyset$.
Assuming the contrary, 
observe that ${\cal H}\cap \Xi$ is $A_0N_0$-invariant
and so $G_0\cdot ({\cal H}\cap \Xi) $ is closed in $\Xi$.
Pick a BHH-function, restrict it
to ${\cal H }\cap \Xi$ and take a minimum point $x_* \in {\cal H}
\cap \Xi $ of the restriction. Then all points of 
the orbit $A_0N_0(x_*)$ are 
minimum points. Therefore $A_0N_0(x_*)$ is
a totally real submanifold of dimension equal to
${\rm dim}\, G_0/K_0 = {\rm dim}_{\Bbb C}\, G/K$ that is contained 
in $\cal H$, contrary
to the fact that ${\cal H}$ is a proper analytic subset.
From these considerations we get the following
description of $\Xi $, see Theorem 8.2.

\medskip
\noindent
{\it Remark.}
For a proof of the inclusion $\Xi \subset \Xi_I$ 
in a more general setting see \cite{M3}. 
Namely, the result is true for a connected real semisimple   
Lie group with two commuting involutions whose product
is a Cartan involution. 
The corresponding fixed point subgroups generalize $G_0$ and $K$. 
The universal domain is defined similarly.
The proof is based on a detailed study of
double coset decompositions. Complex analytic techniques and,
in particular, the existence of BHH-functions are not used.

\medskip
\noindent
{\bf Theorem 8.2.}\ $\Xi = \Xi_I$.

\medskip
\noindent
Since $\Xi_I$ is a connected component of the open
set $\Psi $, which is obtained
by removing a family of
hypersurfaces from the affine variety $G/K$, we see again that $\Xi $ is Stein. 
Since $\Psi $ is the set of all points
for which the $kBk^{-1}$-orbit is open for every $k \in K_0$,
we have
$$\Psi = \{x\in G/K\ \vert \ {\goth g}_x + ({\rm Ad}k)\cdot {\goth b} = \goth g
\ \ {\rm for\ all}\ \ k\in K_0\}.$$
Let $N$ be the normalizer of $K$. Then $\Gamma = N/K$
is a finite group with a free action $x\to x^\gamma $ on $G/K$.  
From the last description of $\Psi $, it follows
that $\Psi ^{\gamma }=
\Psi $ for all $\gamma \in \Gamma $.
Thus $\Gamma $ interchanges the connected components of $\Psi $.
It follows from the definition  
that $\Xi $ 
is contractible, so a non-trivial
finite group cannot act freely on $\Xi$. Hence 
$\Gamma $ interchanges simply transitively the open sets    
$ \Xi ^{\gamma}$.  
Moreover,
for any group $\tilde K \subset G$ with connected component $K$
the covering map $G/K \to G/\tilde K$ induces
a biholomorphic map of $\Xi $ onto its image, cf. \cite{FHW}, Cor. 11.3.6.

\medskip
\noindent
{\bf Theorem 8.3\ {\rm (A.~Huckleberry \cite{Ha}}).}\ 
{\it $\Xi $ is Kobayashi hyperbolic.}

\medskip
\noindent
{\it Proof.}\ By Frobenius reciprocity,
there
exist a $G$-module $V$ and a vector $v_0 \in V$ such that 
$K \subset G_{v_0} \ne G$.
If $G_0$ is of non-Hermitian type then
$K$ is a maximal connected subgroup of $G$.
If $G_0$ is of Hermitian type then there are precisely two intermediate
subgroups between $K$ and $G$, both of them being parabolic. 
In any case,
the connected component of
the stabilizer of the line $[v_0]$ equals $K$
and the natural maps $G/K \to Gv_0 \to G[v_0]$ are finite coverings.
Let ${\Bbb C}[V]_d \subset {\Bbb C}[V]$ 
be the subspace of homogeneous polynomials of degree $d$,
let $I_d$ be the intersection of ${\Bbb C}[V]_d$
with the ideal of (the closure of) $G[v_0]$ 
and let $M_d$ be a $G$-stable complement to $I_d$ in ${\Bbb C}[V]_d$.
The space of all polynomials in $M_d$ vanishing on 
$G[v_0]\setminus B[v_0]$
is $B$-stable and non-trivial for some $d$, so $B$ has 
an eigenvector $\varphi $ in that space.
The zero set of $\varphi $ on the orbit $G[v_0]$
is exactly
the complement to the open $B$-orbit $B[v_0]$.
Replacing $V$ by its symmetric power $S^kV$ and $v_0$ by $v_0^k \in S^kV$,
we obtain a linear form $\varphi $ with the same property.
Now, let $V_0$ be the intersection
of all hypersurfaces $g^* \varphi = 0,\ g\in G$.
Then $V_0$ is a $G$-stable linear subspace of $V$
and we have the $G$-equivariant linear projection map $\pi : V \to W= V/V_0$.
Let $w_0 = \pi (v_0)$ and let $\psi \in W^*$
be the linear form defined by $\pi ^*\psi = \varphi$. Then $K \subset G_{w_0}$
and $G_{w_0} \ne G$,
because $\varphi $ is non-constant on the orbit $Gv_0$. Therefore
$\pi $ gives rise to the finite coverings
$Gv_0 \to Gw_0$ and $G[v_0] \to G[w_0]$. By construction, the orbit
$G\psi =\{g^*\psi \ \vert \ g\in G\}$  
generates
$W^*$ and the same is true for $G_0\psi$. 
By \cite{Ha}, Cor. 2.13, there exist
hyperplanes ${\cal H}_i = \{g_i^*\psi = 0\} \subset {\Bbb P}(W),\
g_i \in G_0,\  i=1, \ldots, 2m + 1,m ={\rm dim}\ 
{\Bbb P}(W)$, satisfying the normal crossing
conditions. It is then known that 
${\Bbb P}(W) \setminus \bigcup_{_i}\, {\cal H}_i$
is Kobayashi hyperbolic, see \cite{K}, Cor. 3.10.9.
The intersection of this set with the orbit $G[w_0]$ equals
$\bigcap _i \, g_iB[w_0]$ and is likewise hyperbolic.
Recall that we have an equivariant fibering $G/K \to G/\tilde K = G[w_0]$.
As we have seen before stating the theorem,
$\Xi $ is mapped biholomorpically onto its image.
The latter is contained in 
the connected component of 
$\bigcap _i \, g_iB[w_0] $ 
at $[w_0]$ and is therefore hyperbolic. 
\hfill $\square $ 

\medskip
\noindent
{\bf Theorem 8.4\ {\rm (G.~Fels, A.~Huckleberry \cite{FH}}).}
{\it If $\Xi ^\prime $ is 
a $G_0$-invariant, Stein,
and Kobayashi hyperbolic domain                                                
in $G/K$ that contains $\Xi $ then
$\Xi ^\prime = \Xi$.}

\medskip
\noindent
The proof requires analysis of the boundary ${\rm bd}(\Xi )$.
First,
one considers 
the special case of $G_0 = SL_2({\Bbb R})$  
and proves Theorem 8.4 for the
crown $\Xi _{{\goth sl}_2}$
of $SL_2({\Bbb R})/SO_2({\Bbb R})$.
Note that $G = SL_2({\Bbb C})$
has precisely two 
non-isomorphic affine homogeneous surfaces.
Namely, if $T \simeq {\Bbb C}^* $
is a maximal torus in $SL_2({\Bbb C})$ 
and $N \subset SL_2({\Bbb C})$ 
is the normalizer of $T$, then these surfaces are
of the form
$Q_1 = {SL_2}({\Bbb C})/T
\simeq ({\Bbb P}^1({\Bbb C})\times 
{\Bbb P}^1({\Bbb C}) ) \setminus \Delta $ and  
$Q_2 = {SL_2}({\Bbb C})/N
 \simeq {\Bbb P}^2({\Bbb C}) \setminus C $,
where $\Delta $ is the diagonal and $C$ is a non-degenerate
curve of degree 2. 
The crown $\Xi _{\goth sl_2}$ 
can be viewed as a domain in $Q_1$ or in $Q_2$.  
In the general case one constructs a $G_0$-stable
open
dense subset ${\rm bd}^{\rm gen}(\Xi ) \subset
{\rm bd}(\Xi)$, such that 
for $z\in {\rm bd}^{\rm gen}(\Xi )$
there exists a simple 3-dimensional subalgebra 
$\goth s _0 \subset \goth g_0$
with the following properties:

\smallskip
(i) the orbit of the corresponding complex group
$S = {\rm exp}(\goth s _0 ^{\Bbb C}) \subset G$ through $z$
is an affine surface, i.e., $Sz \simeq Q_1$ or $Sz \simeq Q_2$;

(ii) under this isomorphism
$Sz \cap \Xi $ is mapped biholomorphically onto $\Xi _{{\goth sl}_2}$. 

\smallskip
\noindent
Now, if $ \Xi^\prime \setminus
\Xi \ne \emptyset$ then one can find a point $z$
as above in $\Xi^\prime \cap {\rm bd}
(\Xi )$. Then $Sz \cap \Xi ^\prime $
properly contains $Sz \cap \Xi$,
contrary to the fact that $\Xi _{\goth sl_2}$ is a maximal 
$SL_2({\Bbb R})$-invariant, Stein
and Kobayashi hyperbolic domain in $Q_1$ or in $Q_2$. 
The details are found in \cite{FHW}, see Thm. 10.6.9.

\medskip
\noindent
{\it Remark.}\ In fact, $\Xi $ is the unique maximal
$G_0$-invariant, Stein, and Kobayashi hyperbolic domain in $G/K$
that contains the base point $o$, see \cite{FHW}, Thm. 11.3.1.
 
\medskip
\noindent
{\it Remark.}\ We refer the reader to \cite{GK} for the
definition of the Shylov-type boundary of the crown
and to \cite{KO} for its simple description and 
applications to
the estimates of automorphic forms. 

\section{The Schubert domain}

We assume here that $G_0$ is of non-Hermitian type. Then the map
$G/K \to G/L$
is a finite covering. 
We have an open $G_0$-orbit $D \subset X = G/P$
and the corresponding compact $K$-orbit $C_0 \subset D$.
Let $q$ denote the complex dimension of $C_0$.
Translations $gC_0, \ g\in G$, are called cycles and
are regarded as points of ${\cal M}_X := G/L$.
The cycle space $\cal{M}_D$
is a domain in ${\cal M}_X$ and the crown $\Xi $
is mapped biholomorphically
onto a domain $\tilde \Xi \subset {\cal M}_X$. We want to prove
the statement of Theorem 7.1, namely, that 
$G\{D\}^\circ $ agrees with $\Omega $. Equivalently, 
we will prove that $\cal{M}_D $ agrees with $\tilde  \Xi  $. 
A.~Huckleberry and J.~Wolf \cite{HW} defined the Schubert domain ${\cal  S}_D$
in ${\cal M}_X$ as follows. Let $B$ be an Iwasawa-Borel subgroup of $G$.
The closures of $B$-orbits on $X$ are called Schubert varieties
(with respect to $B$).
The group $B$ has an open orbit
on any such variety $S$. 
Since the open orbit is affine, 
its complement $S^\prime $ is a hypersurface in $S$.
For topological reasons,
the (finite) set ${\cal S}_{C_0}$ of Schubert varieties
of codimension $q$ intersecting
$C_0$ is non-empty.
One shows that $S^\prime \cap D = \emptyset$
for any $S \in {\cal S}_{C_0}$. Thus the 
incidence variety
$${\cal H}(S):=I(S^\prime )= \{
gC_0 \in {\cal M}_X\ \ \vert \ \
gC_0 \cap S^\prime \ne \emptyset \}$$      
is contained in ${\cal M}_X \setminus {\cal M}_D$.
Clearly, 
${\cal H}(S)$ is $B$-invariant.
Furthermore, one can show
that ${\cal H}(S)$ is an analytic hypersurface in ${\cal M}_X$, see 
\cite{FHW}, Prop.7.4.11. 
For any $k \in K_0$ we have
${\cal M}_D \subset {\cal M}_X \setminus k{\cal H}(S)$.
The set
$$\bigcup _{S \in {\cal S}_{C_0}}\bigl \{ 
\bigcup _{k\in K_0}\, k{\cal H}(S) \bigr \}$$
is closed in ${\cal M}_X$. Its complement is denoted by ${\cal S}_D$
and is called the Schubert domain. By construction,
${\cal S}_D$ is a $G_0$-invariant Stein domain and
$${\cal M}_D \subset {\cal S}_D. \eqno{(*)}$$
On the other hand,   
for any boundary point $z \in {\rm bd}(D)$   
there exist                                            
an Iwasawa decomposition $G_0 = K_0A_0N_0$,
an Iwasawa-Borel subgroup 
$B$ containing $A_0N_0$ and a $B$-invariant variety $S_z$ of
codimension $q + 1 $, such that
$z \in S_z $ and $D \cap S_z = \emptyset$
(a supporting Schubert variety at $z$), see \cite {FHW}, Prop. 9.1.2.
Take a boundary point of
${\cal M}_D$ and consider the
corresponding cycle. It has a point $z \in {\rm bd}(D)$,
hence is contained in the incidence variety
$$I(S_z) := \{gC_0\ \vert \  gC_0 \cap S_z \ne \emptyset\}.$$ 
Obviously, $I(S_z)$ is $B$-invariant and $I(S_z)
\subset {\cal M}_X \setminus {\cal M}_D$,
in particular, $I(S_z) \ne {\cal M}_X$.
But $\tilde \Xi$ is contained
in the open $B$-orbit by Theorem 8.2. Thus
a point of ${\tilde \Xi}$
cannot be a boundary point of ${\cal M}_D$, and it follows that
$$\tilde \Xi \subset {\cal M}_D. \eqno{(**)}$$
Finally, one can modify the proof of Theorem 8.3 to show that
${\cal S}_D$ is hyperbolic. Namely, take the linear
bundle ${\Bbb L}$ over $G/L$ defined by the hypersurface ${\cal H}(S)$,
which appears in the definition of ${\cal S}_D$. 
Then some power ${\Bbb L}^k$ admits a $G$-linearization. Thus
we obtain a non-degenerate equivariant map $G/L
\to
{\Bbb P}(W)$, where a $G$-module $W$ is generated
by a weight vector of $B$. The map is in fact a finite   
covering over the image, which is a $G$-orbit in ${\Bbb P}(W)$
containing the image of ${\cal H}(S)$ as a hyperplane section. 
Since $W$ is irreducible, 
the same argument as in the proof of Theorem 8.3
shows that ${\cal  S}_D$ is hyperbolic.
The inclusions ${(*)}$ and ${(**)}$, together with 
Theorem 8.4, imply
$$\tilde {\Xi} = {\cal M}_D = {\cal S}_D.$$

\section{Complex geometric properties of flag domains}

An open $G_0$-orbit in a complex flag manifold $X = G/P$
is called a flag domain.
One classical example of a flag domain
is a bounded symmetric domain in the dual compact
Hermitian symmetric space.
In this example, a flag domain is a Stein manifold.
However,
this is not the case for an arbitrary flag domain $D$, because
$D$ may contain compact complex submanifolds of
positive dimension.
As we have seen, 
the cycle space of $D$ is always Stein.
Here, we consider the properties of $D$ itself. 

\medskip
\noindent
An open orbit $D =G_0(x_0) \subset X$ is said to be measurable if $D$ carries
a $G_0$-invariant volume element. 
We retain the notation of Section 2.
In particular, $x_0 = e\cdot P,\,\goth p = \goth p_{\Phi }$,
where $\goth p$ and $\Phi $ satisfy (i), (ii) of Theorem 2.4. 

\medskip
\noindent
{\bf Theorem 10.1\ {(\rm J.~Wolf \cite{W}, Thm. 6.3)}.}
{\it The open orbit
$G_0(x_0)$ is measurable if and only if
$\tau \Phi ^r = \Phi ^r $ and $\tau \Phi ^u = - \Phi ^u$.
Equivalently, $G_0(x_0)$ 
is measurable if and only if $\goth p \cap \tau\goth p $
is reductive.}

\medskip
\noindent
Since two fundamental Cartan subalgebras in $\goth g_0$
are conjugate by an inner automorphism of $G_0$,
it follows from the above condition and from
Theorem 2.4 that all open $G_0$-orbits on $X$
are measurable or non-measurable simultaneously.
The proof of Theorem 10.1 can be also found in \cite{FHW}, Sect. 4.5.

\medskip
\noindent
{\bf Example 10.2.}\ Let $\goth g_0$ be a real form of inner type.
Since the Cartan subalgebra $\goth t \subset \goth g$ contains 
a compact Cartan subalgebra $\goth j \subset \goth g_0$, it follows that
$\tau (\alpha ) = -\alpha $ for any root $\alpha $. Thus
the open orbit $G_0(x_0)$
is measurable.

\medskip
\noindent
{\bf Example 10.3.}\ If $P = B$ is a Borel subgroup of $G$, 
then $\Phi = \emptyset $,
$\Phi ^u = \Sigma ^+$ and $\tau \Phi ^u = - \Phi ^u$. Therefore an open
$G_0$-orbit in $G/B$ is measurable. 

\medskip
\noindent
A complex manifold $M$ is said to be $q$-complete
if there is a smooth non-negative
exhaustion function $\varrho : M \to {\Bbb R}$,
whose Levi form has at least $n - q$ positive eigenvalues
at every point of $M$. A fundamental theorem
of A.~Andreotti and H.~Grauert says
that for any coherent sheaf ${\cal F}$ on a $q$-complete
manifold and for all $k > q$
one has
$H^k(M,{\cal F})= 0$, see \cite{AnG}.
Note that in 
the older literature including \cite{AnG}
the manifolds that we call $q$-complete
were called
$(q+1)$-complete.

\medskip
\noindent
{\bf Theorem 10.4 \ {(\rm W.~Schmid, J.~Wolf \cite{SW})}.}
{\it If $D$ is a measurable open $G_0$-orbit in a flag manifold of $G$ and
$q$ is the dimension of the compact $K$-orbit in $D$ then
$D$ is $q$-complete. In particular, $H^k(D, {\cal F}) = 0$
for all coherent sheaves on $D$ and for all $k > q$.}

\medskip
\noindent
The authors of \cite{SW} do not say that $D$ is measurable, but they
use the equivalent condition that the isotropy group of $D$
is the centralizer of a torus.  
The proof of Theorem 10.4 can be also found in \cite{FHW},
see Thm. 4.7.8.

\medskip
\noindent
{\bf Example 10.5.}\ Let $X = {\Bbb P}^n({\Bbb C}),\
G =  SL_{n+1}({\Bbb C})$, and $G_0 =  SL_{n+1}({\Bbb R})$.
Let $\{e_1, e_2, \ldots, e_{n+1}\}$ be a basis of ${\Bbb R}^{n+1}$.
If $n > 1$ then $G_0$ has two orbits on $X$, the open one
and the closed one, with representatives
$x_0 = [e_1+ie_2]$ and $[e_1]$, respectively. 
The isotropy subgroup $(G_0)_{x_0}$ is not reductive. Its 
unipotent radical consists of all $g\in G_0$, such that 
$$g(e_i) = e_i\ (i=1,2),\  g(e_j) \equiv  e_j\ {\rm mod }\,({\Bbb R}e_1 + 
{\Bbb R}e_2)\ (j\ge 3).$$ 
Hence the open orbit $D=G_0(x_0)={\Bbb P}^n({\Bbb C}) \setminus
{\Bbb P}^n({\Bbb R})$
is not measurable. Note that $K =  SO_{n+1}(\Bbb C)$.
Thus the
compact $K$-orbit $C_0 \subset D$ is the projective
quadric $z_1^2 + z_2^2 + \ldots +
z_{n+1}^2 = 0$ and its dimension equals $n-1$. 
In this case, we have $n - q = 1$ and we show how to construct
a smooth non-negative
exhaustion function $\varrho: {\Bbb P}^n({\Bbb C})\setminus 
{\Bbb P}^n({\Bbb R})
\to {\Bbb R}$, whose Levi form has at least one positive eigenvalue
at every point.
For $z=x+iy \in {\Bbb C}^{n+1}$ put
$$\varrho _1(z) = \sum\, x_k^2 + \sum\, y_k^2,\ \ \varrho_2(z) =
\sqrt {\sum\, (x_ky_l -x_ly_k)^2} $$ 
and
notice that
$$\varrho _1(\zeta z) = \vert \zeta \vert^2 \varrho _1(z),\ \ 
\varrho _2(\zeta z) = \vert \zeta \vert ^2 \varrho _2 (z) \ \ {\rm for\ any}
\ \zeta \in {\Bbb C}^*.$$
Thus 
$$\varrho ([z]) = {\varrho _1(z) \over \varrho _2(z)}$$
is well-defined for all $[z] \in {\Bbb P}^n({\Bbb C}) 
\setminus {\Bbb P}^n({\Bbb R})$.
Obviously, $\varrho $ is a smooth exhaustion function for ${\Bbb P}^n({\Bbb C})
\setminus {\Bbb P}^n{(\Bbb R})$.
Given a point $[z] = [x+iy] \in {\Bbb P}^n({\Bbb C})\setminus
{\Bbb P}^n({\Bbb R})$,
take the line $L$ in ${\Bbb P}^n({\Bbb C})$,
connecting $[z]$ with $[x] \in {\Bbb P}^n({\Bbb R})$,
and restrict $\varrho $ to that line.
Clearly, $L$ is the projective image of the affine line
$$\lambda = \alpha + i\beta \mapsto w = x+i\lambda y = x - \beta y 
+ i\alpha y$$
and the restriction $\varrho \vert _L$ equals 
$$\varphi (\lambda) := \varrho ([w]) = {\vert \alpha \vert \over D}\sum y_k^2
\ + {1\over \vert \alpha \vert D}\sum (x_k-\beta y_k)^2,$$
where $D = \varrho _2(x+iy)$.
Computing the Laplacian 
$$\Delta \varphi = {\partial ^2 \varphi \over \partial \alpha ^2} + 
{\partial ^2 \varphi \over \partial \beta ^2}$$
for $\alpha \ne 0 $, we get
$$\Delta \varphi = {2\over D\vert \alpha \vert ^3}
\sum \, (x_k-\beta y_k)^2
+ {2\over D\vert\alpha \vert} \sum \, y_k^2 \, > 0,$$
showing that $\varphi $ is strictly subharmonic.
Hence the Levi form of $\varrho $
has at least one positive eigenvalue at $[z] = [w]\vert_{\alpha =1, \beta =0}$.

\medskip
\noindent 
{\it Concluding remark.} The open orbit in the last
example is not measurable. As a matter of fact, the conclusion
of Theorem 10.4 holds true in this
case. In general, the author does not know whether one can drop the
measurability assumption in Theorem 10.4.

\bigskip
\bigskip
\noindent
{\sl D.N.Akhiezer}

\noindent
Institute for Information Transmission Problems,

\noindent
19 B.Karetny per., Moscow, 127994, Russia.

\noindent

\noindent
Email: akhiezer@iitp.ru

\end{document}